\documentstyle[fleqn,12pt]{article}
\begin{document}
\pagenumbering{arabic}
\setcounter{page}{1}
\pagestyle{plain}
\baselineskip=16pt

\thispagestyle{empty}
\vspace{1.4cm}

\begin{center}
{\large\bf Two-parameter differential calculus on the $h$-superplane}
\end{center}

\vspace{0.5cm}
\noindent
Salih Celik \footnote{E-mail: scelik@fened.msu.edu.tr}\\
\footnote{New E-mail: sacelik@yildiz.edu.tr}
{\footnotesize Mimar Sinan University, Department of Mathematics, 
80690 Besiktas, Istanbul, TURKEY.}

\noindent
Sultan A. Celik \\
{\footnotesize Yildiz Technical University, Department of Mathematics, 
Sisli, Istanbul, TURKEY. }

\vspace{2cm}
{\bf Abstract}

We introduce a noncommutative differential calculus on the two-parameter 
$h$-superplane via a contraction of the $(p,q)$-superplane. 
We manifestly show that the differential calculus is covariant under 
$GL_{h_1,h_2}(1\vert 1)$ transformations. We also give a two-parameter 
deformation of the $(1+1)$-dimensional phase space algebra. 

\vfill\eject
\noindent
{\bf I. INTRODUCTION}

Quantum groups are a generalization of the concept of classical groups. The 
theory of quantum groups has become an important branch of mathematical 
physics and a new branch of mathematics. An approach to obtain the quantum 
groups is to identify the elements of a quantum group with the linear 
transformations of a space with noncommuting coordinates. It is known, from 
the work of Woronowicz,$^1$ that one can define a consistent differential 
calculus on the noncommutative space of a quantum group. Thus 
quantum group is a concrete example of noncommutative differential 
geometry.$^2$ 

During the past few years, Wess-Zumino$^3$ have developed a differential 
calculus on the quantum (hyper)plane which is covariant under the 
action of the quantum group $GL_q(n)$. The natural extension of their 
scheme to superspace$^4$ was given by Soni$^5$ and the two-parameter 
differential calculus on the superplane has been worked out by Chung.$^6$ 
A differential calculus on the $h$-plane was given by Karimipour$^7$ and 
the two-parameter analogue was introduced by Aghamohammadi.$^8$ 

In this paper we construct a two-parameter differential calculus on the 
quantum $h$-superplane using the methods of Ref. 9. The paper is organized 
as follows: in Sec. II we obtain the $(h_1,h_2)$-superplanes via a contraction 
from the $(p,q)$-superplanes. We define derivatives and differentials on the 
$(h_1,h_2)$-superplane of noncommuting coordinates and give their 
commutation rules. In Sec. III we manifestly show that the differential 
calculus is covariant under the action of the quantum supergroup 
$GL_{h_1,h_2}(1\vert 1)$ of Ref. 10. 
We give a two-parameter deformation of the $(1+1)$-dimensional phase 
space algebra in Sec. IV and in the following section we show that the 
$(p,q)$-deformed superoscillator algebra satisfies the undeformed 
superoscillator algebra when objects are transformed into new objects 
such that they are singular for certain values of the deformation 
parameters. 

\noindent
{\bf II. DIFFERENTIAL CALCULUS ON $h$-SUPERPLANE} 

In this work we denote $(p,q)$-deformed objects by primed quantities. 
Unprimed quantities represent transformed coordinates. As usual, 
we shall always assume that even (bosonic) objects commute with everything 
and odd (Grassmann) objects anticommute among themselves. Before discussing
the two-parameter differential calculus on $h$-superplane we give some 
notations and useful formulas in the following section. This first section 
closely follows the approach of Ref. 10. 

\noindent
{\bf A. Quantum $h$-superplane }

Quantum superplane is an associative coordinate algebra ${\cal A}_q$ equipped 
with a set $\{x', \theta'\}$ of generators $x'$, $\theta'$. 
The commutation relations of the generators is defined by$^4$ 
$$ x' \theta' - q \theta' x' = 0, \eqno(1\mbox{a}) $$
$$\theta'^2 = 0, \eqno(1\mbox{b})$$
where $q$ is a nonzero complex deformation parameter. The coordinates 
neither commute nor anticommute unless $q \longrightarrow \pm 1$, 
respectively. In this work we shall use the limits $p \longrightarrow 1$, 
$q \longrightarrow 1$ to make a contraction. 

We now introduce new coordinates $x$ and $\theta$, in terms of $x'$ 
and $\theta'$, by$^{10}$
$$x = x' - {{h_1}\over{p - 1}} \theta', $$
$$  \theta = - {{h_2}\over{q - 1}} x' + 
  \left(1 - {{h_1 h_2}\over{(p - 1)(q - 1)}}\right) \theta'. \eqno(2)$$
Using relation (1), it is easy to verify that 
$$x \theta = q \theta x + h_2 x^2, \eqno(3\mbox{a})$$
where the new deformation parameter $h_2$ commutes with the coordinate $x$ 
and anticommutes with the coordinate $\theta$. Similarly, from (1a) one 
obtains 
$$\theta^2 = - h_2 \theta x, \eqno(3\mbox{b})$$
where 
$$h_1 h_2 = - h_2 h_1, \qquad h_1^2 = 0 = h_2^2. \eqno(4)$$
That is, the new deformation parameters $h_1$ and $h_2$ are odd (Grassmann) 
numbers which anticommute. Note that although in the $p \longrightarrow 1$, 
$q \longrightarrow 1$ limits the transformation (2) is ill behaved, the 
resulting commutation relations are well defined. 

The relations (3) define a new deformation,$^{11}$ which we called the 
$h_2$-deformation of the algebra of coordinate functions on the Manin 
superplane generated by $x$ and $\theta$ in the limit $q \longrightarrow 1$, 
and will be denoted by ${\cal A}_{h_2}$. 

Differential calculus on the quantum superplane ${\cal A}_{h_2}$ requires 
the introduction of differentials ${\sf d}x$, ${\sf d}\theta$. 
The complete framework also includes the commutation relations of these 
differentials with the coordinates and derivatives. 

\noindent
{\bf B. Relations between coordinates and differentials} 

To establish a noncommutative differential calculus on the quantum 
superplane ${\cal A}_{h_2}$, we assume that the commutation relations 
between the coordinates and their differentials have the following form: 
$$x' {\sf d}x' = A {\sf d}x' x', $$
$$x' {\sf d}\theta' = C_{11} {\sf d}\theta' x' + C_{12} {\sf d}x' \theta', 
  \eqno(5)$$
$$\theta' {\sf d}x' = C_{21} {\sf d}x' \theta' + C_{22} {\sf d}\theta' x', $$
$$\theta' {\sf d}\theta' = B {\sf d}\theta' \theta'. $$
Now we would like to transform these relations to unprimed quantities 
to determine the coefficents $A$, $B$ and $C_{ij}$. We first introduce the 
exterior differential {\sf d}. 

The exterior differential {\sf d} is an operator which gives the mapping from 
the coordinates to the differentials 
$${\sf d}: Z^i \longrightarrow {\sf d}Z^i, \eqno(6)$$
where $Z^1 = x$, $Z^2 = \theta$, ${\sf d}Z^1 = {\sf d}x$ and 
${\sf d}Z^2 = {\sf d}\theta$. We demand that the exterior differential {\sf d} 
has to satisfy two properties: the nilpotency 
$${\sf d}^2 = 0, \eqno(7)$$
and the graded Leibniz rule 
$${\sf d}(F G) = ({\sf d} F) G + (- 1)^{\hat{F}} F ({\sf d} G), \eqno(8)$$
where $\hat{F}$ is the Grassmann degree of $F$, that is, 
$\hat{F} = 0$ for even variables and $\hat{F} = 1$ for odd variables. 
We wish to substitute into (5) the differentials ${\sf d}x'$ and 
${\sf d}\theta'$ together with the coordinates $x'$ and $\theta'$. 
The deformation parameters $h_1$ and $h_2$ are both odd numbers and the 
exterior differential {\sf d} is also odd. Therefore the action of the 
exterior differential {\sf d} on $\alpha u$ is defined by 
$${\sf d}(\alpha u) = (- 1)^{\hat{\alpha}} \alpha {\sf d}u, \eqno(9)$$
where $\alpha$ is a number (even or odd) and $u$ is a coordinate of 
superplane. So we can write from (2) 
$$ {\sf d}x = {\sf d}x' + {h_1\over {p - 1}} {\sf d}\theta', $$ 
$${\sf d}\theta = {{h_2}\over {q - 1}} {\sf d}x' + \left(1 - 
  {{h_1h_2}\over {(p - 1)(q - 1)}}\right) {\sf d}\theta'. \eqno(10)$$
Note that if we consider $x$ and $\theta$ as functions of two variables 
(say $x'$ and $\theta'$) and differentiate (2), as usual, then we do 
not obtain the expressions in (10). To obtain (10) one must take 
the differential from the left in Eq. (2). In the Appendix, we 
explain this in detail. 

We now substitute (2) and (10) into (5) which are not explicity written here. 
It will be calculate the coefficients $A$, $B$ and $C_{ij}$. 
We first assume that 
$${\sf d}x' {\sf d}\theta' = p^{- 1} {\sf d}\theta' {\sf d}x', \qquad 
  ({\sf d}x')^2 = 0. \eqno(11)$$
Then we have 
$${\sf d}\theta {\sf d}x = p {\sf d}x {\sf d}\theta - h_1 ({\sf d}\theta)^2, 
  \eqno(12\mbox{a})$$
and 
$$ ({\sf d}x)^2 = h_1 {\sf d}x {\sf d}\theta. \eqno(12\mbox{b})$$
Consequently, the coefficients are determined as follows: 
$$A ~~\mbox{undetermined}, \qquad B = 1, $$
$$C_{11} = q, \qquad C_{12} = pq - 1, \qquad C_{21} = 0, 
  \qquad C_{22} = - p.   \eqno(13)$$
Here we shall choose $A$ equal to $pq$ since the relations are then well 
defined. 

\noindent
{\bf C. Relations of derivatives and coordinates} 

In this section we shall define the derivatives and find the commutation 
relations of derivatives with coordinates and the commutation relations 
between derivatives. We first introduce the matrix$^{10}$ 
$$g = \left(\matrix{ 1 + h_1 h_2/(p-1)(q-1) & h_1/(p-1) \cr 
                     h_2/(q-1) & 1 \cr }\right). \eqno(14)$$
It is easy to verify that the matrix $g$ is a supermatrix. Thus we can write 
the transformation in (2) of the form 
$$Z' = g Z, \qquad Z' = \left(\matrix{x' \cr \theta' \cr}\right). 
  \eqno(15)$$
Let us denote the partial derivatives with respect to $x'$ and $\theta'$ by 
$$\partial_{x'} = {\partial\over {\partial x'}}, \qquad 
  \partial_{\theta'} = {\partial\over {\partial \theta'}}, \eqno(19)$$
respectively. The transformation law of the partial derivatives is then 
defined by 
$$\partial' = (g^{st})^{- 1} \partial, \qquad 
  \partial = \left(\matrix{\partial_x \cr \partial_\theta\cr}\right) 
   \eqno(16)$$
where $g^{st}$ denotes the supertranspose of $g$. Explicitly 
$$\partial_{x'} = \partial_x - {{h_2}\over {q - 1}} \partial_\theta, \quad 
  \partial_{\theta'} = {{h_1}\over {p - 1}} \partial_x + 
  \left(1 - {{h_1 h_2}\over {(p - 1)(q - 1)}}\right) \partial_\theta. 
  \eqno(17)$$
Note that, when one demands the validity of the chain rule, to obtain the 
expressions in (17) it must be assumed that the derivatives act from the 
left on the transformed variables. This case will also be explained in 
detail in the Appendix. 

We know that the exterior differential {\sf d} is defined by 
$${\sf d} = {\sf d}x' \partial_{x'} + {\sf d}\theta' \partial_{\theta'}. 
  \eqno(18\mbox{a})$$
Substituting (10) and (17) into (18a) one obtains 
$${\sf d} = {\sf d}x \partial_x + {\sf d}\theta \partial_\theta, 
  \eqno(18\mbox{b})$$
that is, {\sf d} preserves its form. So, since 
$${\sf d} F(x,\theta) = {\sf d}x \partial_x F + 
  {\sf d}\theta \partial_\theta F, \eqno(19)$$
for any function $F$, replacing $F$ with $x F$ and $\theta F$ we get the 
following relations: 
$$ \partial_x x = 1 + pq x \partial_x + h_1 \theta \partial_x + 
   h_2 x \partial_\theta + 
   h_1 h_2 (x \partial_x + \theta \partial_\theta) + 
   (pq - 1) \theta \partial_\theta, $$
$$\partial_x \theta = p \theta \partial_x - 
  p h_2 (x \partial_x + \theta \partial_\theta),  \eqno(20)$$
$$\partial_\theta x = q x \partial_\theta - q h_1 (x \partial_x + 
      \theta \partial_\theta), $$
$$\partial_\theta \theta =  1 - \theta \partial_\theta + 
   h_1 \theta \partial_x + h_2 x \partial_\theta + 
   h_1 h_2 (x \partial_x + \theta \partial_\theta).$$

We now find the commutation rules between derivatives. These rules can be 
easily obtained by using the nilpotency of the exterior differential. Thus 
we write 
$$0 = {\sf d}^2 = {\sf d}x {\sf d}\theta (p \partial_x \partial_\theta - 
  \partial_\theta \partial_x + h_1 \partial_x^2) + 
  ({\sf d}\theta)^2 (\partial_\theta^2 - h_1 \partial_x \partial_\theta) $$
which says that 
$$\partial_\theta \partial_x = p \partial_x \partial_\theta + 
  h_1 \partial_x^2, \qquad 
  \partial_\theta^2 = h_1 \partial_x \partial_\theta. \eqno(21)$$

The complete framework of the differential calculus requires commutation 
relations of the differentials with derivatives. 

\noindent
{\bf D. Relations of differentials with derivatives}

Finally we shall find the commutation relations between differentials 
and derivatives. We assume that they have the following form in terms of 
primed quantities: 
$$ \partial_{x'} {\sf d}x' = A_{11} {\sf d}x' \partial_{x'} + 
   A_{12} {\sf d}\theta' \partial_{\theta'}, $$
$$ \partial_{x'} {\sf d}\theta' = A_{21} {\sf d}\theta' \partial_{x'} + 
   A_{22} {\sf d}x' \partial_{\theta'}, \eqno(22)$$
$$ \partial_{\theta'} {\sf d}x' = B_{11} {\sf d}x' \partial_{\theta'} + 
   B_{12} {\sf d}\theta' \partial_{x'}, $$
$$\partial_{\theta'} {\sf d}\theta' = B_{21} {\sf d}\theta' \partial_{\theta'} 
  + B_{22} {\sf d}x' \partial_{x'}. $$
Substituting (10) and (17) into (22) and using 
$${\sf d}~ ({\sf d}x) = - ({\sf d}x)~ {\sf d}, \qquad 
  {\sf d}~ ({\sf d}\theta) = ({\sf d}\theta)~ {\sf d}, \eqno(23\mbox{a})$$
and the relation 
$$\partial_i (X^j {\sf d}X^k) = \delta^i{}_j \delta^k{}_l {\sf d}X^k, 
  \eqno(23\mbox{b})$$
where $\partial_1 = \partial_x$ and $\partial_2 = \partial_\theta$, 
we determine the coefficients $A_{ij}$ and $B_{ij}$. So one has 
$$ \partial_{x} {\sf d}x = pq {\sf d}x \partial_x + 
   h_1 {\sf d}\theta \partial_x - h_2 {\sf d}x \partial_\theta + 
   h_1 h_2 ({\sf d}x \partial_x + {\sf d}\theta \partial_\theta) + 
   (pq - 1) {\sf d}\theta \partial_\theta, $$
$$ \partial_x {\sf d}\theta = p {\sf d}\theta \partial_x + 
   p h_2 \left({\sf d}x \partial_x + {\sf d}\theta \partial_\theta\right), $$
$$ \partial_\theta {\sf d}x = - q {\sf d}x \partial_\theta - 
   q h_1 \left({\sf d}x \partial_x + {\sf d}\theta \partial_\theta\right), 
   \eqno(24)$$
$$ \partial_\theta {\sf d}\theta = {\sf d}\theta \partial_{\theta} - 
   h_1 {\sf d}\theta \partial_x + h_2 {\sf d}x \partial_\theta + 
   h_1 h_2 ({\sf d}x \partial_x + {\sf d}\theta). $$

\noindent
{\bf E. Algebra of one-forms}

In this section we shall define two one-forms using the generators of $\cal A$ 
and find the commutation relations of one-forms. 

If we call them $w$ and $u$ then one can define them as follows: 
$$w = {\sf d} x ~x^{-1}, \qquad 
  u = {\sf d} \theta ~x^{-1} - {\sf d} x^{-1} \theta x^{-1}. \eqno(25)$$
We denote the algebra of one-forms generated by two elements $w$ and $u$ by 
$\Omega$. The generators of the algebra $\Omega$ with the generators of 
$\cal A$ satisfy the following relations: 
$$x w = w x - h_1 u x, \qquad \theta w = - w \theta + h_1 u \theta, $$
$$x u = u x, \qquad \theta u = u \theta - h_2 (w \theta + u x). \eqno(26)$$

The commutation rules of the generators of $\Omega$ are 
$$w^2 = 0, \qquad w u = u w. \eqno(27)$$

Using (18b) and (25), if we define the operators $T$ and $\nabla$ as 
$$T = x \partial_x + \theta \partial_\theta, 
  \qquad \nabla = x \partial_\theta, \eqno(28)$$
then we have 
$$T \nabla = \nabla T, \qquad \nabla^2 = 0, \eqno(29)$$
as a subalgebra of gl$(1\vert 1)$. 

The action of $T$ and $\nabla$ on the generators $x$ and $\theta$ is 
$$T x = x + x T, \qquad \nabla x = x \nabla - h_1 x T, $$
$$T \theta = \theta + \theta T, \qquad 
  \nabla \theta = x - \theta \nabla + h_1 \theta T. \eqno(30)$$

\noindent\vfill\eject
{\bf III. THE SUPERGROUP $GL_{h_1,h_2}(1\vert 1)$ AND COVARIANCE} 

It is well known that the quantum supergroup $GL_{p,q}(1\vert 1)$ acts 
as a linear transformation on the quantum superplane, preserves (1) 
and the dual relations 
$$\varphi'^2 = 0, \qquad \varphi' y' - p^{-1} y' \varphi' = 0. \eqno(31)$$
In extending this property of covariance under the coaction of 
$GL_{p,q}(1\vert 1)$, from the superplane to its calculus, it will be assumed 
that the deformed group structure implies and is implied by invariance 
of the intermediary relations (5) under linear transformations of the 
quantum superplane. In the present work, this will be applied to the 
$(h_1,h_2)$-deformed superplane.  

In this section we would like to discuss the meaning of covariance in a 
graded version of noncommutative differential calculus of Wess-Zumino$^3$ 
for the two-parameter case. Before proceeding, we define the dual quantum 
$h$-superplane. 

To define the dual quantum superplane, we interpret the 
differentials ${\sf d}x$ and ${\sf d}\theta$, as the coordinates of the dual 
superplane, as follows 
$${\sf d}x = \varphi, \qquad {\sf d}\theta = y. \eqno(32)$$
Now the quantum dual $h$-superplane generated by $y$, $\varphi$ with the 
relations (12) in the limit $p \longrightarrow 1$ will be denoted by 
${\sf d} {\cal A}_{h_1}$. If we assume that ${\cal A}_{h_2}$ and 
${\sf d} {\cal A}_{h_1}$ have to be covariant under the coaction 
$$\delta(x) = a \otimes x + \beta \otimes \theta, \qquad 
  \delta(\theta) = \gamma \otimes x + d \otimes \theta, \eqno(33\mbox{a})$$ 
$$\delta({\sf d} Z) = (\tau \otimes {\sf d}) \delta(Z), \quad 
  \tau(u) = (-1)^{\hat u} u \eqno(33\mbox{b})$$ 
and that $\beta$, $\gamma$ anti-commute with $\theta$, $\varphi$, $h_1$ and 
$h_2$ we get the corresponding $(h_1,h_2)$-deformation of the supergroup 
$GL(1\vert 1)$ as a quantum matrix supergroup $GL_{h_1,h_2}(1\vert 1)$ 
generated by $a$, $\beta$, $\gamma$, $d$ with the relations$^{10}$ 
$$a \beta = \beta a - h_1 (a^2 - \beta \gamma - ad), \quad 
  d \beta = \beta d + h_1 (d^2 + \beta \gamma - da), $$
$$a \gamma = \gamma a + h_2 (a^2 + \gamma \beta - ad), \quad 
  d \gamma = \gamma d - h_2 (d^2 - \gamma \beta - da), $$
$$ \beta^2 = h_1 \beta (a - d), \quad \gamma^2 = h_2 \gamma (d - a), $$
$$ \beta \gamma = - \gamma \beta + (h_1 \gamma - h_2 \beta)(a - d), 
   \eqno(34)$$
$$ ad = da + h_1 (a - d) \gamma + h_2 \beta (a - d), $$
where 
$${\cal D} = ad^{-1} - \beta d^{-1} \gamma d^{-1} = 
d^{-1} a - d^{-1} \beta d^{-1} \gamma. $$

The two-parameter differential calculus on the quantum superplane is 
explicitly as follows: 

The commutation relations of variables and their differentials are 
$$x \theta = \theta x + h_2 x^2, \qquad \theta^2 = - h_2 \theta x, $$
$$\varphi y = y \varphi + h_1 y^2, \qquad \varphi^2 = h_1 \varphi y. 
  \eqno(35)$$
Note that the last two relations of (35) are obtained from (14) and 
(15). Hovewer they can also be obtained from (31) with the limits 
$p \longrightarrow 1$, $q \longrightarrow 1$. 

The commutation relations between variables and derivatives are 
$$\partial_x x = 1 + x \partial_x - h_1 \theta \partial_x +
  h_2 x \partial_\theta + h_1 h_2 (x \partial_x + \theta \partial_\theta), $$
$$\partial_x \theta = \theta \partial_x - 
  h_2 (x \partial_x + \theta \partial_\theta), \eqno(36) $$
$$\partial_\theta x = x \partial_\theta - 
  h_1 (x \partial_x + \theta \partial_\theta), $$
$$\partial_\theta \theta = 1 - \theta \partial_\theta - h_1 \theta \partial_x 
  + h_2 x \partial_\theta + h_1 h_2 (x \partial_x + \theta \partial_\theta),$$
and those among the derivatives are 
$$\partial_x \partial_\theta = \partial_\theta \partial_x - 
  h_1 \partial_x^2, \qquad 
  \partial_\theta^2 = h_1 \partial_\theta \partial_x. \eqno(37)$$

The commutation relations of variables with their differentials are 
$$x \varphi = \varphi x + h_1 (\varphi \theta - y x) + h_1 h_2 \varphi x,$$
$$x y = y x - h_1 y \theta - h_2 \varphi x + h_1 h_2 \varphi \theta,$$
$$\theta \varphi = - \varphi \theta + h_1 y \theta - h_2 \varphi x 
  - h_1 h_2 y x,\eqno(38)$$
$$\theta y = y \theta - h_2 (\varphi \theta + y x)- h_1 h_2 y \theta. $$

The commutation relations between derivatives and differentials are 
$$\partial_x \varphi = \varphi \partial_x + h_1 y \partial_x - 
   h_2 \varphi \partial_\theta + 
   h_1 h_2 (\varphi \partial_x + y \partial_\theta),$$ 
$$\partial_x y = y \partial_x + h_2 (\varphi \partial_x + y \partial_\theta), 
  \eqno(39) $$
$$\partial_\theta \varphi = - \varphi \partial_\theta - 
  h_1 (\varphi \partial_x + y \partial_\theta), $$
$$\partial_\theta y = y \partial_\theta - 
  h_1 y \partial_x + h_2 \varphi \partial_\theta +
  h_1 h_2 (\varphi \partial_x + y \partial_\theta). $$
Note that this calculus goes back to those of Ref. 9 when $h_1 = 0$ 
and $h_2 = h$. This calculus is slightly different from Ref. 10. The 
reason for this difference is the use of commutation relations of the dual 
{\it exterior} superplane in Ref. 10 instead of the dual superplane in this 
work. 

We now discuss the covariance of the differential calculus. The covariance 
here means that all the relations between coordinates $x$, 
$\theta$, differentials ${\sf d}x$, ${\sf d}\theta$ and derivatives 
$\partial_x$, $\partial_\theta$, etc. must preserve their form when one 
changes the coordinates by 
$$ x \longrightarrow a x + \beta \theta, \qquad 
 \theta \longrightarrow \gamma x + d \theta, \eqno(40)$$
where the matrix $T = \left(\matrix{a & \beta \cr \gamma & d \cr}\right)$ 
is an element of the quantum supergroup GL$_{h_1,h_2}(1\vert 1)$ 
acting on the quantum superplane. We must change the differentials by 
$${\sf d}x \longrightarrow a {\sf d}x - \beta {\sf d}\theta, \qquad 
  {\sf d}\theta \longrightarrow - \gamma {\sf d}x + d {\sf d}\theta, 
   \eqno(41)$$
since the odd objects anti-commute among themselves. 
Covariance can be maintained if one defines the transformation law of the 
partial derivatives as folllows 
$$\partial_x \longrightarrow (a^{-1} - a^{-1} \gamma d^{-1} \beta a^{-1}) 
  \partial_x - a^{-1} \gamma d^{-1} \partial_\theta, $$
$$\partial_\theta \longrightarrow (d^{-1} - d^{-1} \beta a^{-1} \gamma d^{-1}) 
  \partial_\theta + d^{-1} \beta a^{-1} \partial_x.  \eqno(42) $$

\noindent
{\bf IV. A TWO-PARAMETER DEFORMATION OF CLASSICAL PHASE SPACE}

We shall now give a two-parameter deformation of the 
$(1 + 1)$-dimensional classical phase space. We denote the 
algebra $(35)-(37)$ generated by coordinates $x$, $\theta$ 
and the derivatives $\partial_x$ and $\partial_\theta$ by 
${\cal B}_{h_1,h_2}$. It is interesting to note that simply 
identifying $\partial_x$ and $\partial_\theta$ with $i p_x$ 
and $p_\theta$ is not compatible with the hermiticity of 
coordinates and momenta. To identify $\partial_x$ and $\partial_\theta$ 
with the momenta $i p_x$ and $p_\theta$, one must take care of the 
hermiticity of the coordinates and momenta. To this end, we first define 
the hermitean conjugation of the coordinates $x$ and $\theta$, respectively, 
as 
$$x^+ = (1 + 2h_1 h_2) x + 2h_1 \theta, \quad 
  \theta^+ = (1 - 2h_1 h_2) \theta + 2h_2 x. \eqno(43)$$
It is then easy to see that the hermiticity of $x^+$ and $\theta^+$ impose 
some condition on the deformation parameters, i.e., $h_1$ is a real parameter 
and $h_2$ is a pure imaginary parameter: 
$$\overline{h_1} = h_1, \qquad \overline{h_2} = - h_2,  \eqno(44)$$
where the bar denotes complex conjugation. In this case, the 
hermitean conjugation of the derivatives $\partial_x$ and $\partial_\theta$ 
are 
$$\partial_x^+ = - (1 + 2h_1 h_2) \partial_x + 2h_2 \partial_\theta, \quad 
  \partial_\theta^+ = (1 - 2h_1 h_2) \partial_\theta + 2 h_1 \partial_x. 
  \eqno(45)$$
In the $h_1 \longrightarrow 0$, $h_2 \longrightarrow 0$ limits the 
definitions (43) and (45) go back to those of the classical case. 

The relations $(35)-(37)$ are now invariant under the transformations 
(43) and (45). The above involution allows us to define the hermitean 
operators 
$$\hat{x} = (1 + h_1 h_2) x + h_1 \theta, \quad 
  \hat{\theta} = (1 - h_1 h_2) \theta + h_2 x, \eqno(46)$$
and, as bosonic and fermionic momenta, 
$$\hat{p}_x = i [(1 + h_1 h_2) \partial_x - h_2 \partial_\theta], \quad 
  \hat{p}_\theta = (1 - h_1 h_2) \partial_\theta + h_1 \partial_x. \eqno(47)$$
The final form of the $(h_1,h_2)$-deformed phase space algebra is 
$$\hat{x} \hat{\theta} = \hat{\theta} \hat{x} + h_2 \hat{x}^2, \qquad 
  \hat{\theta}^2 = - h_2 \hat{\theta} \hat{x}, $$
$$\hat{p}_x \hat{p}_\theta = \hat{p}_\theta \hat{p}_x + 
  i h_1 \hat{p}_x^2, \qquad 
   \hat{p}_\theta^2 = - i h_1 \hat{p}_x \hat{p}_\theta, $$
$$\hat{p}_x \hat{x} = i + \hat{x} \hat{p}_x + 
   i h_2 \hat{x} \hat{p}_\theta - h_1 \hat{\theta} \hat{p}_x + 
   h_1 h_2 (1 + \hat{x} \hat{p}_x + i \hat{\theta} \hat{p}_\theta), $$
$$\hat{p}_x \hat{\theta} = \hat{\theta} \hat{p}_x - 
  h_2 (\hat{x} \hat{p}_x + i \hat{\theta} \hat{p}_\theta),  \eqno(48)$$
$$\hat{p}_\theta \hat{x} = \hat{x} \hat{p}_\theta + 
   h_1 (i \hat{x} \hat{p}_x - \hat{\theta} \hat{p}_\theta), $$
$$\hat{p}_\theta \hat{\theta} = 1 - \hat{\theta} \hat{p}_\theta + 
  h_2 \hat{x} \hat{p}_\theta + i h_1 \hat{\theta} \hat{p}_x - 
  h_1 h_2 (1 + i \hat{x} \hat{p}_x - \hat{\theta} \hat{p}_\theta). $$
This gives a $(h_1,h_2)$-deformed phase space algebra which may be used 
to study the $(1+1)$-dimensional quantum phase space. 

Note that we can derive a deformed super-Clifford algebra from the phase 
space algebra as follows: suppose that we define gamma matrices 
$$\gamma^1 \equiv \hat{p}_\theta, \qquad 
  \gamma^2 \equiv \hat{\theta}, \qquad
  c^1 \equiv \hat{p}_x, \qquad c^2 \equiv \hat{x}. \eqno(49)$$
Then, they satisfy super-Clifford algebra 
$$c^1 c^2 = c^2 c^1 - h_1 \gamma^2 c^1 + i(1 + h_2 c^2  \gamma^1) + 
  h_1 h_2 (1 + \gamma^2 \gamma^1 + c^2 c^1), $$
$$c^1 \gamma^2 = \gamma^2 c^1 - h_2 (c^2 c^1 + i \gamma^2 \gamma^1), $$
$$\gamma^1 c^1 = c^1 \gamma^1 - i h_1 (c^1)^2,  \eqno(50) $$
$$\gamma^1 c^2 = c^2 \gamma^1 - h_1 (\gamma^2 \gamma^1 - i c^2 c^1), $$
$$\gamma^1 \gamma^2 = 1 - \gamma^2 \gamma^1 + i h_1 \gamma^2 c^1 + 
   h_2 c^2 \gamma^1 - h_1 h_2 (1 + c^2 c^1 - \gamma^2 \gamma^1), $$
$$(\gamma^1)^2 = - i h_1 c^1 \gamma^1, \qquad 
  (\gamma^2)^2 = - h_2 \gamma^2 c^2, $$
$$\gamma^2 c^2 = c^2  \gamma^2 - h_2 (c^2)^2. $$

\vfill\eject
\noindent
{\bf V. A COMMENT ON SUPEROSCILLATORS} 

We know that introducing one 'bosonic' and one 'fermionic' oscillator, $A$ 
and $B$, respectively, and making the usual identification 
$$x' ~\longleftrightarrow~ A^+, \qquad \theta' ~\longleftrightarrow~ B^+,$$
$$\partial_{x'} ~\longleftrightarrow~ A, \qquad 
  \partial_{\theta'} ~\longleftrightarrow~ B, \eqno(51)$$
one constructs the quantum super-oscillator algebra which is covariant under 
the quantum supergroup GL$_{p,q}(1\vert 1)$. Under identification (2) and 
(17) one has 
$$x ~\longleftrightarrow~ A^+ - {{h_1}\over {p - 1}} B^+, \quad 
  \partial_x ~\longleftrightarrow~ 
  \left(1 + {{h_1 h_2}\over {(p - 1)(q - 1)}}\right) A + 
  {{h_2}\over {q - 1}} B, $$
$$ \theta ~\longleftrightarrow~ 
  \left(1 - {{h_1 h_2}\over {(p - 1)(q - 1)}}\right) B^+ 
   - {{h_2}\over {q - 1}} A^+, \quad 
 \partial_\theta ~\longleftrightarrow~ B - {{h_1}\over {p - 1}} A, \eqno(52)$$
where 
$$\overline{p} = q. \eqno(53)$$ 
Substituting (52) into (20) and (3), surprisingly all 
$(h_1,h_2)$-dependence cancels and one obtains the usual 
$(p,q)$-deformed super-oscillator algebra$^{12}$ 
$$A A^+ = 1 + pq A^+ A + (pq - 1) B^+ B, $$
$$B B^+ = 1 - B^+ B, \qquad B^2 = 0 = B^{+ 2}, $$
$$A B^+ = p B^+ A, \qquad A B = p^{- 1} B A, \eqno(54)$$
$$A^+ B = q^{- 1} B A^+, \qquad A^+ B^+ = q B^+ A^+, $$
In the $p \longrightarrow 1$, $q \longrightarrow 1$ limits, we 
get undeformed super-oscillator algebra. 

\noindent
{\bf VI. APPENDIX}

In this Appendix, we show that a two-parameter covariant differential 
calculus on the quantum $h$-superplane can be constructed only if the 
derivatives and differentials act from the left. 

Consider the change of coordinates which is given by (2) 
$$x' = \left(1 + {{h_1 h_2}\over {(p - 1)(q - 1)}}\right) x + 
       {{h_1}\over {p - 1}} \theta, \quad 
  \theta' = \theta + {{h_2}\over {q - 1}} x. \eqno(\mbox{A}1)$$
If we interpret the symbols ${\sf d}x$ and ${\sf d}\theta$ as 
differentials acting from the right and demanding the validity of the 
chain rule, we have 
$${\sf d}x = {\sf d}x' {{\partial x}\over {\partial x'}} + 
  {\sf d}\theta' {{\partial x}\over {\partial \theta'}} = 
{\sf d}x' + {\sf d}\theta' \left(- {{h_1}\over {p - 1}}\right) = 
 {\sf d}x' - {{h_1}\over {p - 1}}{\sf d}\theta' \eqno(\mbox{A}2)$$
and 
$${\sf d}\theta = {{h_2}\over {q - 1}} {\sf d}x' + 
 \left(1 + {{h_1 h_2}\over {(p - 1)(q - 1)}}\right) {\sf d}\theta'. 
 \eqno(\mbox{A}3)$$
Therefore, for example, 
$${\sf d}x' = \left(1 + {{h_1 h_2}\over {(p - 1)(q - 1)}}\right) {\sf d}x + 
  {{h_1}\over {p - 1}} {\sf d}\theta  \eqno(\mbox{A}4)$$
so that 
$$\mbox{RHS of (A4)} ~ = 
  \left(1 + {{h_1 h_2}\over {(p - 1)(q - 1)}}\right) {\sf d}x' - 
  {{h_1}\over {p - 1}} {\sf d}\theta' + 
   {{h_1}\over {p - 1}} {\sf d}\theta' + 
   {{h_1 h_2}\over {(p - 1)(q - 1)}} {\sf d}x' \neq {\sf d}x'. $$

Similarly, if we write, from the chain rule, 
$$\partial_x = 
  \left(1 + {{h_1 h_2}\over {(p - 1)(q - 1)}}\right) \partial_{x'} + 
  {{h_2}\over {q - 1}} \partial_{\theta'}, \eqno(\mbox{A}5)$$
and 
$$\partial_\theta = {{h_1}\over {p - 1}} \partial_{x'} + \partial_{\theta'}, 
  \eqno(\mbox{A}6)$$
then 
$$\partial_{x'} = \partial_x - {{h_2}\over {q - 1}} \partial_\theta, \quad
  \partial_{\theta'} = - {{h_1}\over {p - 1}} \partial_x + 
  \left(1 - {{h_1 h_2}\over {(p - 1)(q - 1)}}\right) \partial_\theta, 
\eqno(\mbox{A}7)$$
so that, for example 
$$\mbox{RHS of (A6)} ~ = \left(1 - 2{{h_1 h_2}\over {(p - 1)(q - 1)}}\right) 
  \partial_\theta \neq \partial_\theta. $$

This asymmetry between right and left derivative and differential for 
transformed variables stems from the matrix $g$ in (14) which off 
diagonal elements are odd. That is, $g$ is a supermatrix and so the 
supertranspose must be used. 

\noindent
{\bf ACKNOWLEDGMENT}

This work was supported in part by T.B.T.A.K the Turkish 
Scientific and Technical Research Council. 

\vfill\eject
\noindent
{\footnotesize 
$^1$ S. L. Woronowicz, Comm. Math. Phys. {\bf 122}, 125 (1989). \\
$^2$ A. Connes, A., {\it Non-commutative differential geometry}, 
     (Academic, London, 1994). \\
$^3$ J. Wess and B. Zumino, Nucl. Phys. Proc. Suppl. B {\bf 18}, 302 (1990). \\
$^4$ Yu I. Manin, Comm. Math. Phys. {\bf 123}, 163 (1989). \\
$^5$ S. Soni, J. Phys. A {\bf 24}, L459 (1991). \\
$^6$ W. S. Chung, J. Math. Phys. {\bf 35}, 2484 (1994). \\
$^7$ V. Karimipour, Lett. Math. Phys. {\bf 30}, 87 (1994). \\
$^8$ A. Aghamohammadi, Mod. Phys. Lett. A {\bf 8}, 2607 (1993). \\
$^9$ S. Celik, S. A. Celik, and M. Arik, 
     J. Math. Phys. {\bf 39}, 3426 (1998). \\
$^{10}$ S. Celik, Lett. Math. Phys. {\bf 42}, 299 (1997). \\
$^{11}$ L. Dabrowski and P. Parashar, Lett. Math. Phys. {\bf 38}, 331 (1996).\\
$^{12}$ S. Vokos, J. Math. Phys. {\bf 32}, 2979 (1991). }

\end{document}